\documentclass[a4paper]{amsart}
\usepackage{amsmath}
\usepackage{amsfonts}
\usepackage{amssymb}
\usepackage{cite}

\newtheorem{theorem}{Theorem}[section]

\newtheorem{claim}[theorem]{Claim}
\newtheorem{corollary}[theorem]{Corollary}
\newtheorem{proposition}[theorem]{Proposition}

\theoremstyle{definition}

\newtheorem{notation}[theorem]{Notation}

\DeclareMathOperator\tp{tp}

\newcommand\mycoloneqq{\mathrel{\mathop:}=}

\begin{document}
\title{Uniform bounds on growth in o-minimal structures}
\author{Janak~Ramakrishnan} 

\begin{abstract}
We prove that a function definable with parameters in an o-minimal structure is
bounded away from $\infty$ as its argument goes to $\infty$ by a function
definable without parameters, and that this new function can be chosen
independently of the parameters in the original function.  This generalizes a
result in \cite{FM5}.  Moreover, this remains true if the argument is taken to
approach any element of the structure (or $\pm\infty$), and the function has
limit any element of the structure (or $\pm\infty$).

\end{abstract}
\maketitle

\section{Introduction}

We begin with a special case of the main result of this paper.

\begin{proposition}\label{propparam} Let $M$ be an o-minimal expansion of a
  dense linear order $(M,<)$.  Let $f:M^n\times M\longmapsto M$ be definable in
  $M$.  Then there exist functions $g:M\longmapsto M$ and $h:M^n\longmapsto M$
  definable in $M$ such that $f(x,t)\le g(t)$ for all $x\in M^n$ and $t>h(x)$.
  Moreover, if $M'$ is the prime model containing the parameters used to define
  $f$, then $g$ and $h$ are defined over $M'$.  \end{proposition}

This was already known under the additional assumption that $M$ expands an
ordered group; see 3.1 of \cite{FM5}, which uses \cite[C.4]{vdDMi96} and
\cite{MiSt98}.  Here, we remove the need for the group structure.  Indeed, we
show something stronger.

\begin{theorem}\label{paramsbounded} Let $M$ be an
  o-minimal expansion of a dense linear order $(M,<)$.  Let $f$ be an $n+1$-ary
  $M$-definable function with domain $A\times M$ for some $A\subseteq M^n$.
  Suppose that, for some $b\in M\cup\{\infty\}$ and all $x\in A$, we have
  $\lim_{t\to b^-}f(x,t)=b$ and $f(x,t)<b$.  Then there exist functions
  $g:M\longmapsto M$ and $h:A\longmapsto M$ definable in $M$ such that $h(A)<b$
  and for $t\in(h(x),b)$, we have $g(t)\in [f(x,t),b)$.  Moreover, if $M'$ is
  the prime model containing the parameters used to define $f$, then $g$ and $h$
  are defined over $M'$.\end{theorem}

If $M$ expands a field, then using the maps $1/(b-t)$ and $b-1/t$ this theorem
follows easily from 3.1 of \cite{FM5}.  When $b=\infty$ and $M$ expands an
ordered group, this is essentially 3.1 of \cite{FM5}.  However, this result is
new if $M$ does not expand a group, or if $M$ does not expand a field and $b\in
M$.

Corollary \ref{cor} strengthens the theorem slightly, allowing $f$ to take any
value as its limit, from either direction.  Note that if Corollary \ref{cor} is
applied in the case that $f$ is definable in the prime model of an o-minimal
theory, this shows that any definable function is bounded as it approaches a
limit by one definable in the prime model, assuming the limit is in the prime
model or $\pm\infty$.

We use the terminology of \cite{Tressl05a}: the definable $1$-types in an
o-minimal theory are called ``principal.''  To each principal type over a
structure $M$ is associated a unique element $a\in M\cup\{\pm\infty\}$ to which
it is ``closest,'' in the sense that no elements of $M$ lie between $a$ and any
realization of the type.  We say that a principal type is ``principal
above/below/near $a$.''  We write $\langle a_1,\ldots,a_n\rangle$ to denote the
tuple of length $n$ having the element $a_i$ as its $i$th component.

\section{Results}

\begin{proof}[Proof of Theorem \ref{paramsbounded}]
We first note that the theorem is equivalent to the following:

\begin{claim}\label{modelclaim} Let $P$ be the prime model of the theory of $M$,
  let $b\in P\cup\{\infty\}$, and let $A\subseteq M^n$ be a
  $\emptyset$-definable set.  Let $f(x,t):A\times M\longmapsto M$ be a
  $\emptyset$-definable function and $a\in A$ a tuple, with $\lim_{t\to
    b^-}f(a,t)=b$ and $f(a,t)<b$ for all $t<b$.  Then there exists a
  $\emptyset$-definable function $g:M\longmapsto M$ such that $g(t)\in [f(a,t),b)$ for
  $t$ sufficiently close to $b$.  Similarly if the limit is taken as $t$
  approaches $b$ from above and $f$ approaches $b$ from above, with $b\in
  M\cup\{-\infty\}$.  \end{claim}

The theorem implies Claim \ref{modelclaim}, since the fact that $g$ bounds $f$
transfers to elementary extensions.  Inversely, if the theorem failed, then we
could add the parameters needed to define $f$ to the language, so that $f$ would
become $\emptyset$-definable, and then by compactness we could find $a$ in an
elementary extension of $P$ such that the claim failed.  (Note that $b$ is
definable from the parameters used to define $f$.)  Therefore, we prove Claim
\ref{modelclaim}.

\begin{notation}
  Let $w$ be an element in some elementary extension of $M$ that realizes the
  principal type below $b$ over $M$.  In other words, $w$ is infinitesimally
  close to $b$ with respect to $M$.  It is easy to see that all $M$-definable
  functions extend to this elementary extension, and that if $\varphi$ is any
  $\emptyset$-definable (respectively $M$-definable) predicate, $\varphi(w)$
  holds if and only if $\varphi(t)$ holds for all $t$ in some interval $(c,b)$,
  with $c\in P$ (respectively $c\in M$).  Thus, whenever we write $\varphi(w)$,
  the reader should understand this as equivalent to ``$\varphi(t)$ for all $t$
  in some interval with right endpoint $b$ and left endpoint definable over the
  same parameters used to define $\varphi$.''  \end{notation}

We go by induction on the length of $a$, simultaneously for all o-minimal
structures, all $\emptyset$-definable functions, and all tuples of appropriate
length.  Let $f(x,t)$ and $a$ satisfy the conditions of Claim \ref{modelclaim}
for some $b$.  If $a=\langle a_1,\ldots,a_n\rangle$ with $n>1$, we can add
constants for $a_1,\ldots,a_{n-1}$ to the language and use induction for the
cases of $n-1$ and $1$ to prove the claim.  Thus, we may suppose that $a$ is a
singleton.  If $a\in P$, then the claim is trivial, so suppose not.

We can use regular cell decomposition \cite[2.19(2)]{vdDries98} to ensure that
$f$ is monotone in $x$ and increasing in $t$ on its two-dimensional domain cell,
$C$, which we can take to be
\begin{equation*}
\{\langle x,t\rangle\mid x\in(d_1,d_2)\land k(x)<t<b\},
\end{equation*}
for some $\emptyset$-definable monotone function $k$ and $d_1,d_2\in
P\cup\{\pm\infty\}$ (with $d_1<a<d_2$).  We may also require that $f(C)<b$.

The case where $f(x,w)$ is constant in $x$ at $a$ is easy by standard
o-minimality arguments, since then the value $f(a,w)$ is definable from $w$
without using $a$.  Thus, we may suppose that $f(x,t)$ is non-constant in $x$ at
$a$, for all $t\in(k(a),b)$.  Without loss of generality, assume that $f$ is
increasing in $x$ on $C$.

If $\tp(a)$ is not principal below $d_2$, then we can choose $a'\in P$ with
$a<a'<d_2$.  Then $f(a',t)> f(a,t)$ for $t\in(\max\{k(a),k(a')\},b)$, and so we
are done.  Thus, we may suppose that $\tp(a)$ is principal below $d_2$.

The proof relies on the following claim.

\begin{claim}\label{tpinter}
Let $p\in S_1(\emptyset)$ be the principal type below $b$.  If there is no
$\emptyset$-definable map between $\tp(a)$ and $p$, then Claim \ref{modelclaim}
holds.
\end{claim}
\begin{proof}
  If $k(a)\models p$, then $k$ is the desired map between $\tp(a)$ and $p$.
  Thus, we can assume that $k(a)<c$ for some $c\in P$ with $c<b$.  Increasing
  $d_1$ if necessary, we may also assume that $c\ge\sup\{k(x)\mid x\in
  (d_1,d_2)\}$, so if $t\in (c,b)$ and $x\in(d_1,d_2)$, then $\langle
  x,t\rangle\in C$.  Now consider the formula
\begin{equation*}
    \varphi(t)\mycoloneqq \sup\{f(x,t)\mid x\in(d_1,d_2)\}=b.
\end{equation*}
First, suppose that $\varphi(w)$ does not hold.  Then, for any $t$
sufficiently close to $b$,
\begin{equation*}
  \sup\{f(x,t)\mid x\in (d_1,d_2)\}<b.
\end{equation*}
Let $z(t)$ be this (uniformly $t$-definable) supremum.  Then
$z(t)\in[f(a,t),b)$, and so Claim \ref{modelclaim} holds.  Thus, the case that
remains to consider is when $\varphi(w)$ does hold.  We can then fix
$t_0\in(c,b)$ with $t_0\in P$ such that $\varphi(t_0)$ holds, and we have a
$\emptyset$-definable map, $f(x,t_0)$.  We show $f(a,t_0)\models p$.  For any
$e\in P$ with $e<b$, we can find $r\in (d_1,d_2)\cap P$ such that
$f(r,t_0)\in(e,b)$, by $\varphi(t_0)$.  Since $r<a$ (else $a$ would not be
principal below $d_2$) and $f(x,t_0)$ is increasing in $x$, we have
$f(a,t_0)>f(r,t_0)>e$.  Thus, $f(a,t_0)\models p$, witnessing the
$\emptyset$-definable map between $\tp(a)$ and $p$.
\end{proof}

We now complete the proof of Claim 2.1.  By Claim \ref{tpinter}, we can assume
that $\tp(a)$ is principal below $b$.  Then the domain cell $C$ has the form
\begin{equation*} \{\langle x,t\rangle\mid x\in(d_1,b)\land k(x)<t<b\}.
\end{equation*}

If $k(w)\ge f(a,w)$, then we are done, so we may assume that $f(a,w)>k(w)$.
Then we may increase $d_1$ and suppose that for any $x\in(d_1,b)$ we have
$f(x,w)>k(w)$, as well as $f(x,w)>w$.  Fix $e\in P$ with $e\in(d_1,b)$.  We have
$f(e,w)>k(w)$.  Then $\langle w,f(e,w)\rangle\in C$.  For any $t\in(a,b)$, since
$f(e,t)>t$ and $f$ is increasing in both coordinates, $f(t,f(e,t))>f(a,t)$.  So
we are done, since $f(t,f(e,t))$ is $\emptyset$-definable and
$f(t,f(e,t))\in(f(a,t),b)$ for $t$ sufficiently close to $b$ -- namely, for
$t\in(a,b)$.  \end{proof}

\begin{corollary}\label{cor} Theorem \ref{paramsbounded} holds when $\lim_{t\to
    b^{\pm}}f(x,t)=c$, with $c\in M\cup\{\pm\infty\}$ and $f$ approaching $c$
  from either direction, with $g$ and $h$ now definable over the prime
  model containing $c$ and the parameters defining $f$.  \end{corollary}
\begin{proof} Suppose that the limit is taken as $t$ approaches $b$ from below
  and that $f(x,t)$ approaches $c$ from above.  The other cases are similar.
  Let $P$ be the prime model of the statement.  Choose $a\in A\cap P^n$.  Let
  $\psi(t)$ denote the inverse of $f(a,t)$, so $\psi$ is $\emptyset$-definable.
  Then for any $x\in M^n$, the limit of $\psi(f(x,t))$ is $b$ as $t$ goes to
  $b$, and this value approaches $b$ from below.  By Theorem
  \ref{paramsbounded}, there are $\emptyset$-definable $\tilde g$ and $\tilde h$
  with $\tilde g(t)\in [\psi(f(x,t)),b)$ for $t\in(\tilde h(x),b)$.  Then for
  $t\in(\tilde h(x),b)$, we have $f(a,\tilde g(t))\in(c,f(x,t)]$.  Since
  $f(a,\tilde g(t))$ is still $\emptyset$-definable, $f(a,\tilde g(t))$ is the
  desired function $g$, and $\tilde h$ is the desired function $h$.  \end{proof}

\bibliography{../janak}
\bibliographystyle{halpha}

\end{document}